\documentclass[11pt,leqno]{article}%
\usepackage{amsmath}
\usepackage{amsthm}
\usepackage{amssymb}
\usepackage{amsfonts}
\usepackage{graphicx}%
\providecommand{\U}[1]{\protect\rule{.1in}{.1in}}
\newtheorem{thm}{Theorem}[section]
\newtheorem{lma}{Lemma}[section]
\theoremstyle{definition}
\newtheorem*{defn}{Definition}
\newtheorem*{example}{Example}
\newtheorem*{ack}{Acknowledgment}
\theoremstyle{remark}
\newtheorem*{rmk}{Remark}
\newenvironment{pf}{\begin{proof}}{\end{proof}}
\makeatletter\@addtoreset{equation}{section}\makeatother

\begin{document}

\title{Rate of Decay of Stable Periodic Solutions of Duffing Equations}
\author{Hongbin Chen \thanks{Supported in part by the science funds of Xi'an Jiaotong
University }\\Department of Mathematics\\Xian Jiaotong University\\Xian, P.R.~China 029-82667909
\and Yi Li \thanks{Supported in part by the Xiao-Xiang Funds, Hunan Normal
University }\\Department of Mathematics \\University of Iowa \\Iowa City, Iowa 52242, USA  \\and \\Department of Mathematics \\Hunan Normal University \\Changsha, Hunan, China\\ }
\date{}
\maketitle

\begin{abstract}
In this paper, we consider the second-order equations of Duffing type. Bounds
for the derivative of the restoring force are given that ensure the existence
and uniqueness of a periodic solution. Furthermore, the stability of the
unique periodic solution is analyzed; the sharp rate of exponential decay is
determined for a solution that is near to the unique periodic solution.

\end{abstract}

\begin{flushleft}
Key words: Periodic solution; Stability; Rate of decay.

AMS subject classification: 34B15, 34D05.
\end{flushleft}

\section{\label{int}Introduction and statement of main results}

This paper is devoted to the existence, uniqueness and stability of periodic
solutions of the Duffing-type equation
\begin{equation}
\label{i1}x^{\prime\prime}+cx^{\prime}+g(t,x)=h(t),
\end{equation}
where $g\colon\mathbb{R}\times\mathbb{R}\rightarrow\mathbb{R}$ is a
$T$-periodic function in $t$ and continuous in $x$, $h(t)$ is a $T$-periodic
function, and $c >0$ is a positive constant. The existence and multiplicity of
periodic solutions of (\ref{i1}) or more general types of nonlinear
second-order differential equations have been investigated extensively by many
authors. However, the stability of periodic solutions is less extensively
studied. In \cite{Ort}, R.~Ortega studied (\ref{i1}) from the stability point
of view and obtained stability results by topological index. A.C.~Lazer and
P.J.~McKenna get stability results by converting the equation (\ref{i1}) to a
fixed-point problem \cite{LaMc}. Recently, more complete results concerning
the stability of periodic solutions of (\ref{i1}) were obtained by J.M.~Alonso
and R.~Ortega \cite{Al1,Al2}. Under the condition that the derivative of the
restoring force is independent of $t$ and positive, they found sharp bounds
that guarantee global asymptotic stability. In \cite{Al2}, optimal bounds for
stability are obtained. However, the rate of exponential decay of periodic
solutions is not known in general. The purpose of this paper is to determine
the rate of decay of the periodic solutions. Our results are motivated by
careful observation of the following example of a second-order differential
equation with constant coefficients:
\begin{equation}
\label{i2}x^{\prime\prime}+cx^{\prime}+kx=h(t),
\end{equation}
where $c>0$. If $k\neq0$, then (\ref{i2}) has a unique $T$-periodic solution
$x_{0}(t)$, and all solutions of (\ref{i2}) have the form $x(t)=c_{1}%
e^{\rho_{1}t}+c_{2}e^{\rho_{2}t}+x_{0}(t)$ where $\rho_{1}<\rho_{2}$ are roots
of $\rho^{2}+c\rho+k=0$. Regarding $k$ as a parameter, when $k<0$, then
$\rho_{2}>0$, and the unique periodic solution is unstable; when $k>0$, the
unique periodic solution is stable; keep increasing $k$ till it crosses the
critical value $k=\frac{c^{2}}{4}$, then $\rho_{1}$, $\rho_{2}$ are a pair of
conjugates. In this case, every solution other than the unique $T$-periodic
solution decays to the periodic solution at the same exponential rate
$\frac{c}{2}$ independently of $k$.

The aim of this paper is to show that the above-mentioned $\frac{c}{2}$ decay
results can be generalized to (\ref{i1}). In order to show that periodic
solutions of (\ref{i1}) are stable, it is essential to impose conditions on
the restoring force that can rule out the existence of additional periodic
solutions that are subharmonic of order $2$, whereas to determine the rate of
decay of periodic solutions, it is essential to impose conditions on $g(x,t)$
such that the corresponding Dirichlet boundary-value problem does not admit
any nontrivial solutions. We shall give existence and uniqueness results,
characterizing solutions that are locally asymptotically stable with sharp
exponential rate of decay, by a Sturm comparison argument coupled with Floquet
theory. Roughly speaking, (\ref{i1}) has a unique $T$-periodic solution which
is locally asymptotically stable in the sense of Lyapunov, as long as the
restoring force is small or the fractional constant is large.\medskip

The following notations will be used throughout the rest of the paper.\medskip

1.\quad$L_{T }^{p}$ \qquad$T$-periodic function $u\in L^{p}[\,0,T \,]$ with
$\Vert u\Vert_{p}$ for $1\leq p\leq\infty$;

2.\quad$C_{T }^{k}$ \qquad$T$-periodic function $u\in C^{k}[\,0,T \,]$,
$k\geq0$, with $C^{k}$-norm;

3.\quad$\alpha(t)\gg\beta(t)$, if $\alpha(t)\geq\beta(t)$ and $\alpha
(t)>\beta(t)$ on some positive-measure subset.\medskip

Now we state our main results.

\medskip

\begin{thm}
\label{thmb} Assume that $g(t,x)\in C^{1}(\mathbb{R})$ satisfies the following conditions:

\begin{enumerate}
\item \label{thmb(1)} $g^{\prime}_{x}(t,x)<\frac{\pi^{2}}{T^{2}}+\frac{c^{2}%
}{4}$, and

\item \label{thmb(2)} there is an $\alpha(t)\in C_{T}$ with $\overline
{\alpha(t)}\gg\frac{c^{2}}{4}$ such that $g^{\prime}_{x}(t,x)\gg\alpha(t)$,
where $\overline{\alpha(t)}$ denotes the the average of $\alpha(t)$ over a period.
\end{enumerate}

Then \textup{(\ref{i1})} has a unique $T$-periodic solution which is
asymptotically stable with sharp rate of decay of $\frac{c}{2}$ for $c>0$.
\end{thm}

Here, we say that the periodic solution $u_{0}$ of (\ref{i1}) is locally
asymptotically stable if there exist constants $C>0$ and $\alpha>0$ such that
if $u$ is another solution with $|u(0)-u_{0}(0)|+|u^{\prime}(0)-u_{0}^{\prime
}(0)|=d$ sufficiently small, then $|u(t)-u_{0}(t)|+|u^{\prime}(t)-u_{0}%
^{\prime}(t)|<Cde^{-\alpha t}$. The super exponent $\alpha$ as above is called
the rate of decay of $u_{0}$. \medskip

\begin{rmk}
The bounds given in Theorem \ref{thmb} are optimal. In fact, the example given
at the beginning of the section shows that the lower bound is sharp. An
example to show that the upper bound is optimal in the theorem will be given
in the following.
\end{rmk}

\begin{example}
Consider the linear differential equation
\begin{equation}
x^{\prime\prime}+cx^{\prime}+\frac{1}{4}(1+c^{2}+\epsilon\cos t)x=0,
\label{i4}%
\end{equation}
where $c>0$ and $|\epsilon|$ is small and $\epsilon\neq0$. By the
transformation $y(t)=e^{-\frac{c}{2}t}x(t)$, the damping term can be
eliminated and equation \textup{(\ref{i4})} is equivalent to
\begin{equation}
y^{\prime\prime}+\frac{1}{4}(1+\epsilon\cos t)y=0. \label{i5}%
\end{equation}
The following notations are used. Obviously, the nontrivial solution $x(t)$ of
\textup{(\ref{i4})} has rate of decay $\frac{c}{2}$ if and only if the
corresponding solution $y(t)$ of \textup{(\ref{i5})} is a nontrivial bounded
solution. Let $y_{i}$, $i=1,2$, denote the solutions of \textup{(\ref{i5})}
with initial values
\[
y_{1}(0,\epsilon)=y_{2}^{\prime}(0,\epsilon)=1,\qquad y_{1}^{\prime
}(0,\epsilon)=y_{2}(0,\epsilon)=0.
\]

If the discriminant function of \textup{(\ref{i4})} is denoted by
$\triangle(\epsilon)=y_{1}(2\pi,\epsilon)+y_{2}^{\prime}(2\pi,\epsilon)$, then
the Floquet multipliers are the roots of the quadratic
\[
\mu^{2}-\triangle(\epsilon)\mu+1=0.
\]

By means of perturbation one can compute that
\[
\triangle(\epsilon)=-2-\frac{\pi}{64}\epsilon^{2}+0(\epsilon^{2})
\]
for $\epsilon$ small, which shows that the modulus of one of the multipliers
is greater than $1$ and the modulus of the other multiplier is less than $1$.
Therefore the multipliers of \textup{(\ref{i5})} are a pair of distinct real
numbers. This implies that the rate of decay of \textup{(\ref{i5})} is greater
than $\frac{c}{2}$. Thus if the derivative of the restoring term crosses over
the bound given in Theorem \textup{\ref{thmb}}, the conclusion of Theorem
\textup{\ref{thmb}} does not hold.
\end{example}

Consider the piecewise linear equation
\begin{equation}
\label{ii}x^{\prime\prime}+cx^{\prime}+a(t)x^{+}-b(t)x^{-}=h(t),
\end{equation}
where $x^{+}=\max\{x,0\}$, $x^{-}=\min\{-x,0\}$, $a(t), b(t)\in C_{T}$. For
$a$ and $b$ constant, (\ref{ii}) is a simple mathematical model for vertical
oscillations of a long-span suspension bridge, which has received great
attention after a series of works of Lazer and McKenna \cite{LaMck}.

\begin{thm}
\label{thm:a} Suppose that $h(t)\in C_{T}$ has a finite number of zeros in
$[\,0,T\,]$. Then

\begin{enumerate}
\raggedright

\item \label{thm:a(1)} equation \textup{(\ref{ii})} has a unique $T$-periodic
solution if $0\ll a(t),b(t)\ll\frac{(2\pi)^{2}}{T^{2}}+\frac{c^{2}}{4}$,

\item \label{thm:a(2)} the unique $T$-periodic solution is asymptotically
stable if $0\ll a(t),b(t)\ll\frac{(\pi)^{2}}{T^{2}}+\frac{c^{2}}{4}$, and

\item \label{thm:a(3)} the unique $T$-periodic solution has rate of decay
$\frac{c}{2}$ if $\frac{c^{2}}{4}\ll a(t),b(t)\ll\frac{(\pi)^{2}}{T^{2}}%
+\frac{c^{2}}{4}$.
\end{enumerate}
\end{thm}

\section{\label{LinPer}The linear periodic problem}

In this section we shall recall some basic results about topological methods
and prove some stability results for linear periodic systems.

Consider the periodic boundary-value problem
\begin{equation}
\label{hbaa}\left\{
\begin{array}
[c]{l}%
x^{\prime}=F(t,x),\\
x(0)=x(T).
\end{array}
\right.
\end{equation}
where $F\colon[\,0,T \,]\times\mathbb{R}^{n}\rightarrow\mathbb{R}^{n}$ is a
continuous function that is $T$-periodic in $t$. In order to use a homotopic
method to compute the degree, we assume that $h\colon[\,0,T\,]\times
\mathbb{R}^{n}\times[\,0,1\,]\rightarrow\mathbb{R}^{n}$ is a continuous
function such that
\begin{align*}
h(t,x,1)  &  =F(t,x),\\
h(t,x,0)  &  =G(x),
\end{align*}
where $G(x)$ is continuous. The following theorem is due to J.~Mawhin
\cite{Maw}.

\begin{lma}
\label{lmaa} Let $\Omega\subset C_{T }$ be an open bounded set such that the
following conditions are satisfied.

\textup{(1)} There is no $x\in\partial\Omega$ such that
\[
x^{\prime}=h(t,x,\lambda)\;\;\;\;\;\forall\, \lambda\in[\,0,1).
\]

\textup{(2)} $\deg(g,\Omega\cap\mathbb{R}^{n},0)\neq0$.

\noindent Then \textup{(\ref{hbaa})} has at least one solution.
\end{lma}

Next we consider the above system for $n=2$. We denote by $x(t,x_{0})$ the
initial-value solution of (\ref{hbaa}) and introduce the Poincar\'{e} map
$P\colon x_{0}\rightarrow x(T, x_{0})$. It is well known that $x(t,x_{0})$ is
a $T$-periodic solution of the system (\ref{hbaa}) if and only if $x_{0}$ is a
fixed point of $P$. If $x$ is an isolated $T$-periodic solution of
(\ref{hbaa}), then $x_{0}$ is an isolated fixed point of $P$. Hence the
Brouwer index is defined by
\[
\operatorname{ind}[\,P,x_{0}\,]=\deg(I-P,B_{\varepsilon}(0),0).
\]

\begin{defn}
A $T$-periodic solution $x$ of (\ref{hbaa}) will be called a nondegenerate
$T$-periodic solution if the linearized equation
\begin{equation}
\label{hbab}y^{\prime}=f_{x}(t,x)y
\end{equation}
does not admit any nontrivial $T$-periodic solutions.
\end{defn}

Let $M(t)$ be the fundamental matrix of (\ref{hbab}) and $\mu_{1} $ and
$\mu_{2}$ the eigenvalues of the matrix $M(T )$. Then $x(t,x_{0})$ is
asymptotically stable if and only if $\vert\mu_{i}\vert<1$, $i=1,2$:
otherwise, if there is an eigenvalue of $M(T )$ with modulus greater than one,
then $x(t,x_{0})$ is unstable. \medskip

Before giving some results concerning stability of linear periodic equations,
we consider the following eigenvalue problem:
\begin{equation}
\label{hbac}\left\{
\begin{array}
[c]{l}%
x^{\prime\prime}+\lambda x=0,\\
x(0)=x(T)=0,\quad\operatorname{sgn}x^{\prime}(0)= \operatorname{sgn}x^{\prime
}(T).
\end{array}
\right.
\end{equation}

It is easy to verify that $\lambda_{n}=\frac{(2n\pi)^{2}}{T^{2}}$, the $n$-th
eigenvalue of (\ref{hbac}) with the eigenfunction $\varphi_{n}(t)=\sin
\frac{2n\pi}{T}t$. Now for
\begin{equation}
\left\{
\begin{array}
[c]{l}%
x^{\prime\prime}+q(t)x=0,\\
x(0)=x(T)=0,\quad\operatorname{sgn}x^{\prime}(0)=\operatorname{sgn}x^{\prime
}(T),
\end{array}
\right.  \label{hbad}%
\end{equation}
we have the following result.\medskip

\begin{lma}
\label{lmabB}Assume that $q(t)\ll\frac{(2\pi)^{2}}{T^{2}}$. Then equation
\textup{(\ref{hbad})} does not admit any nontrivial solutions.
\end{lma}

\begin{pf}
Suppose that $x(t)$ is a nontrivial solution of (\ref{hbad}). Then $x$ is an
eigenfunction. There is an $n$ such that $\lambda_{n}(q)=0$. Since
$q(t)\ll\frac{(2\pi)^{2}}{T^{2}}$, by a theorem concerning comparison of
eigenvalues, we have that $\lambda_{n}(q)>\lambda_{n}\left(  \frac{(2\pi)^{2}%
}{T^{2}}\right)  =\lambda_{n}-\lambda_{1}>0$, a contradiction.
\end{pf}

Consider the homogeneous periodic equation
\begin{equation}
\label{hbae}L_{\alpha}x:=x^{\prime\prime}+cx^{\prime}+\alpha(t)x=0,
\end{equation}
where $c\in\mathbb{R} $ is constant and $\alpha(t)\in L_{T }^{\infty}$.

\medskip

\begin{lma}
\label{lmac} Assume $\alpha(t)\in L_{T }^{\infty}$ satisfies the following
conditions: $\alpha(t)\ll\frac{(2\pi)^{2}}{T^{2}}+\frac{c^{2}}{4}$ and
$\overline{\alpha(t)}>0$.

Then \textup{(\ref{hbae})} does not admit any nontrivial $T$-periodic solutions.
\end{lma}

\begin{pf}
Suppose on the contrary that (\ref{hbae}) admits a nontrivial $T$-periodic
solution $x(t)$. We claim that $x(t)$ vanishes at some $t_{0}\in[\,0,T\,]$. If
not, then $x(t)\neq0$ for all $t$ in $\mathbb{R}$. By the periodic boundary
conditions, we have $x^{\prime}(T)= x^{\prime}(0)$ and $\frac{x^{\prime}%
(T)}{x(T)}=\frac{x^{\prime}(0)}{x(0)}$. Dividing (\ref{hbae}) by $x(t) $ and
integrating by parts gives that
\[
\int_{0}^{T}\frac{x^{\prime}(t)^{2}}{x(t)^{2}}\,dt +\int_{0}^{T}%
\alpha(t)\,dt=0,
\]
which contradicts the hypothesis of the lemma. So $x(t)$ has a zero in
$[\,0,T\,]$. We may assume that $x(0)=0$ so that $x(0)=x(T)=0$. By the
transformation $y(t)=e^{\frac{c}{2}t}x(t)$, $y(t)$ is a nontrivial solution of
(\ref{hbae}) with $q(t)=\alpha(t)-\frac{c^{2}}{4}$, and if the first condition
of the lemma holds, $q(t)\ll\frac{(2\pi)^{2}}{T^{2}}$. Then according to Lemma
\ref{lmabB}, $y(t)\equiv0$, hence $x(t)\equiv0$, a contradiction.
\end{pf}

The following simple lemma, given in \cite{ChenL}, will be used in proving the
existence and the uniqueness of periodic solutions.

\begin{lma}
\label{lmad}Suppose that $\alpha(t)$, $\alpha_{1}(t)$ and $\alpha_{2}(t)\in
L_{T }^{\infty}$ such that $\alpha_{1}(t)$, $\alpha_{2}(t)$, and $\alpha(t)$
are all $\ll\frac{(2\pi)^{2}}{T^{2}}+\frac{c^{2}}{4}$. Then

\begin{enumerate}
\item \label{lmad(1)} the possible $T$-periodic solution $x$ of equation
\textup{(\ref{hbae})} is either trivial or different from zero for each
$t\in\mathbb{R}$;

\item \label{lmad(2)} $L_{\alpha_{i}}x=0$ \textup{(}$i=1,2$\textup{)} cannot
admit nontrivial $T$-periodic solutions simultaneously if $\alpha_{1}%
(t)\ll\alpha_{2}(t)$; and

\item \label{lmad(3)} $L_{\alpha}x=0$ has no nontrivial $T$-periodic solution,
if $\alpha(t)\gg0$ \textup{(}resp., $\alpha(t)\ll0$\textup{)}.
\end{enumerate}
\end{lma}

The following lemma is essential for determining the rates of decay of
periodic solutions of (\ref{i1}).

\begin{lma}
\label{lmae}Assume that $\alpha(t)\ll\frac{\pi^{2}}{T^{2}}+\frac{c^{2}}{4}$
and that $\overline{\alpha(t)}>\frac{c^{2}}{4}$. Then \textup{(\ref{hbae})}
does not admit real Floquet multipliers.
\end{lma}

\begin{pf}
If the conclusion of the lemma does not hold, then there is a real Floquet
multiplier $\rho$ and a nontrivial solution $x(t)$ such that $x(t+T)=\rho
x(t)$. Introduce the transformation $x=e^{-\frac{1}{2}ct}u$. Then $u$ solves
the equation
\begin{equation}
\label{hbaf}u^{\prime\prime}+\left[  \alpha(t)-\frac{c^{2}}{4}\right]  u=0
\end{equation}
with Floquet multiplier $\chi=e^{-\frac{cT }{2}}\rho$. If $x$ does not change
sign in $\mathbb{R}$, neither does $u$. Dividing equation (\ref{hbaf}) by $u$
and integrating by parts, and noting that $\frac{u^{\prime}(T)}{u(T)}%
=\frac{u^{\prime}(0)}{u(0)}$, we have that
\[
\int_{0}^{T}\frac{u^{\prime}(t)^{2}}{u(t)^{2}}\,dt +\int_{0}^{T}\left[
\alpha(t)-\frac{c^{2}}{4}\right]  \,dt=0,
\]
which contradicts the assumption that $\overline{\alpha(t)}>\frac{c^{2}}{4}$.
Therefore $x(t)$ vanishes at some $t_{0} \in[\,0,T\,]$, and the conditions are
such that $x(t_{0}+T)=\rho x(t_{0})=0$. Thus the corresponding $u$ is a
nontrivial solution of the following Dirichlet boundary-value problem:
\begin{equation}
\label{hbag}u^{\prime\prime}+\left[  \alpha(t)-\frac{c^{2}}{4}\right]  u=0
,u(t_{0})=u(t_{0}+T)=0.
\end{equation}
Since $\alpha(t)-\frac{c^{2}}{4}\ll\frac{\pi^{2}}{T^{2}}$, multiplying
(\ref{hbag}) by $u$ and integrating we have that
\[
\int_{t_{0}}^{T+t_{0}}u^{\prime\,2}\,dt=\int_{t_{0}}^{T+t_{0}}\left[
\alpha(t)-\frac{c^{2}}{4}\right]  u^{2}\,dt <\frac{\pi^{2}}{T^{2}}\int_{t_{0}%
}^{T+t_{0}}u^{2}\,dt,
\]
which contradicts the Poincar\'{e} inequality. Therefore (\ref{hbae}) does not
admit real multipliers.
\end{pf}

\begin{lma}
\label{lmaf} Under the conditions of Lemma \textup{\ref{lmae}}, the rate of
decay of any nontrivial solution of \textup{(\ref{hbae})} is $\frac{c}{2}$.
\end{lma}

\begin{pf}
Consider the corresponding system
\begin{equation}
\label{hbah}X^{\prime}(t)=A(t)X(t),
\end{equation}
where the column vector function $X(t)=(x(t),x^{\prime}(t))^{T}$ and $A(t)$ is
the matrix function
\[
A(t)=\left(
\begin{array}
[c]{cc}%
0 & 1\\
-p(t) & -c
\end{array}
\right)  .
\]
Let $M(t)$ be a fundamental matrix solution of (\ref{hbah}). It is well known
that $M(t)$ has the form
\begin{equation}
\label{hbai}M(t)=P(t)e^{Bt},
\end{equation}
where $P(t)$ and $B$ are $2\times2$ matrices, $P(t)=P(t+T)$, and $B$ is a
constant matrix. Let $\rho_{1}=e^{T\lambda_{1}}$, $\rho_{2}=e^{T\lambda_{2}}$
be the Floquet multipliers and $\lambda_{1}$ and $\lambda_{2}$ the Floquet
exponents associated with $\rho_{1}$ and $\rho_{2}$. Let $x_{1}$ and $x_{2}$
be the eigenvector components of the matrix $e^{TB}$. It follows from Lemma
\ref{lmae} that $\rho_{1}$ and $\rho_{2}$ are a pair of conjugates. Thus the
eigenvectors associated with different eigenvalues are linearly independent.
Therefore $y_{i}=p_{i}(t)e^{\lambda_{i}t}$ (for $i=1,2$) form the fundamental
solution system of equation (\ref{hbah}). On the other hand, by applying the
Jacobi--Liouville formula, we have
\[
\vert\rho_{1}\vert^{2}=\rho_{1}\rho_{2}=e^{-\int_{0}^{T}c\,dt }=e^{-cT}%
\]
and
\begin{align*}
\operatorname{Re}\lambda_{1}  &  =\operatorname{Re}\lambda_{2}=\frac{1}%
{2}\operatorname{Re}\left(  \lambda_{1}+\lambda_{2}\right) \\
&  =\frac{1}{2T}\ln(\rho_{1}\rho_{2}) =-\frac{c}{2}.
\end{align*}
Since every solution is a linear combination of $y_{1}(t)$ and $y_{2}(t)$,
$p_{i}(t)$ is $T$-periodic, hence it is bounded. Therefore every nonzero
solution of equation (\ref{hbae}) decays at the same exponential rate of
$\frac{c}{2}$.
\end{pf}

\section{\label{ProThms}Proof of main results}

Now we prove our main results.

\subsection{\label{Prothma}Proof of Theorem \ref{thmb}}

We begin with the following existence result, dividing the proof into two
steps.\medskip

Step 1. Existence.

Define $F\colon C^{2}_{T}\rightarrow C_{T}$ by
\[
F(x(t)):=x^{\prime\prime}+cx^{\prime}+g(t,x(t)).
\]

We have the following.

\begin{thm}
\label{thmca} Assume that $g(x,t)\in C^{1}(\mathbb{R}\times\mathbb{R})$, and
that $g(x,t)$ is $T$-periodic in $t$. If in addition $g$ satisfies the two conditions

\begin{enumerate}
\item \label{thmca(1)} $g(t,x)/x \ll\frac{(2\pi)^{2}}{T^{2}}+\frac{c^{2}}{4}\;
\forall\, x\in\mathbb{R}$ and

\item \label{thmca(2)} there is a $T$-periodic function $\beta(t)\in C_{T}$
such that $\overline{\beta(t)}>0$ and $g(t,x)/x\gg\beta(t)$ for all
$x\in\mathbb{R}$,
\end{enumerate}

then the differential equation \textup{(\ref{i1})} has a $T$-periodic solution.
\end{thm}

\begin{pf}
Without loss of generality, we may assume that $g(0,t)=0$, for otherwise we
can reduce both sides of equation (\ref{i1}) by $g(0,t)$. Consider the
parametrized equation
\begin{equation}
F_{\lambda}:=x^{\prime\prime}+cx^{\prime}+\lambda g(t,x)+(1-\lambda)ax=\lambda
h(t) \label{hbcb}%
\end{equation}
for some $a\in\left(  0,\frac{(2\pi)^{2}}{T^{2}}+\frac{c^{2}}{4}\right)  $. We
claim that there is an $R>0$ such that equation (\ref{hbcb}) has no solution
on $\partial B_{R}$ for any $\lambda\in\lbrack\,0,1\,]$. If there is not such
an $R$, let $x_{n}$ be a sequence of solutions such that $\Vert x_{n}%
\Vert\rightarrow\infty$ and $\lambda_{n}\in\lbrack\,0,1\,]$. Denote by $z_{n}$
the quotient $z_{n}:=\frac{x_{n}}{\Vert x_{n}\Vert}$. Dividing (\ref{hbcb}) by
$\Vert x_{n}\Vert$, then multiplying by $\varphi(t)\in C_{T}^{2}$ and
integrating by parts, we have that
\begin{equation}
\int_{0}^{T}{z_{n}\varphi^{\prime\prime}-cz_{n}\varphi^{\prime}+g(t,x_{n}%
)\varphi/\Vert x_{n}\Vert}\,dt=\lambda_{n}\int_{0}^{T}\varphi h_{n}/\Vert
x_{n}\Vert\,dt. \label{hbcc}%
\end{equation}

The conditions of Theorem \ref{thmca} imply that $\{[\lambda_{n}%
g(t,x_{n})+(1-\lambda_{n})a]x_{n}/\Vert x_{n}\Vert\}$ is bounded. It is
pre-compact in the weak topology in $L^{1}[\,0,T\,]$. Thus there are
subsequences such that $g(t,x_{n})/x_{n}\rightarrow q(t)$ and $\lambda
_{n}\rightarrow\lambda$. Taking the limit in the equation (\ref{hbcc}), one
obtains that
\begin{equation}
\int_{0}^{T}\{z\varphi^{\prime\prime}-cz\varphi^{\prime}+w(t)\varphi
z\}\,dt=0, \label{hbcd}%
\end{equation}
where $w(t)=\lambda q(t)+(1-\lambda)a$. Evidently, $\overline{w(t)}\geq
\lambda\overline{\alpha(t)}+(1-\lambda)a>0$, and $w(t)<\frac{(2\pi)^{2}}%
{T^{2}}+\frac{c^{2}}{4}$, which satisfies the conditions of Lemma 2.2 in
\cite{ChenL}, and since here $z(t)$ is a $T$-periodic solution, it follows
from Lemma \ref{lmac} that $z(t)\equiv0$, which contradicts the fact that, by
construction, $\Vert z(t)\Vert=1$. Next, by applying the homotopic invariance
property, we have that
\[
\deg(F_{0},B_{R},0)=\deg(F_{1},B_{R},0)=1.
\]
According to Lemma \ref{lmaa}, equation (\ref{i1}) then has a $T$-periodic
solution. This completes the proof of existence.

\medskip

Now, we will finish the proof of Theorem \ref{thmca} by verifying uniqueness.
Let $x_{1}$ and $x_{2}$ be two distinct $T$-periodic solutions of the equation
(\ref{i1}), and let $u=x_{1}- x_{2}$. Then $u$ satisfies the equation
\begin{equation}
\label{hbce}u^{\prime\prime}+cu^{\prime}+p(t)u=0,
\end{equation}
where $p(t)=[g(t,x_{1}(t))-g(t,x_{1}(t))]/[x_{1}- x_{2}]$. If $u\neq0$
identically, it will imply that $u$ is an eigenfunction associated with a
Floquet multiplier equal to one. Again, Lemma \ref{lmac} rules out this
possibility. Therefore $u\equiv0$.
\end{pf}

Step 2. Rate of decay.

To show that every solution of the nonlinear equation (\ref{i1}) locally
decays at the rate of $\frac{c}{2}$ to the unique $T$-periodic solution, we
need the following $C^{1}$ version of the Hartman--Grobman theorem \cite{Bel}.

\begin{lma}
\label{lamcc} Let $f\colon U\subset\mathbb{R}^{n} \rightarrow\mathbb{R}^{n}$
be a $C^{1}$ function with $f(0)=0$ such that $f^{\prime}_{x}(0)\colon
\mathbb{R}^{n}\rightarrow\mathbb{R}^{n}$ is a contracting mapping. Then $f$ is
$C^{1}$ conjugate equivalent to $f^{\prime}_{x}(0)$.
\end{lma}

Consider the planar system associated with equation (\ref{i1}),
\begin{equation}
\label{hbcj}\left\{
\begin{array}
[c]{rl}%
x^{\prime} & =cx-y,\\
y^{\prime} & =h(t)-g(t,x).
\end{array}
\right.
\end{equation}

Let $X_{0}(t)=(x_{0}(t),y_{0}(t))$ be the unique $T$-periodic solution
determined by the initial condition $X_{0}(0)=(x_{0},y_{0})$. Then $X_{0}$
corresponds to the unique fixed point of the Poincar\'{e} mapping $PX=U(T,X)$,
where $U(t,X)$ is the initial-value solution of (\ref{hbcj}) with $U(0,X)=X $.
Let $M(t) $ be the fundamental matrix solution of the linearization
\begin{equation}
\label{hbck}X^{\prime}=A(t)X
\end{equation}
of (\ref{hbcj}), where
\[
A(t)=\left(
\begin{array}
[c]{cc}%
c & -1\\
-q(t) & 0
\end{array}
\right)  .
\]
By the differentiability of $X(t)$ with respect to the initial value, the
Poincar\'{e} mapping can be expressed in terms of the initial value $X$ by the
following formula:
\begin{equation}
\label{hbcl}PX-X_{0}=M(T)(X-X_{0})+o(X-X_{0}).
\end{equation}

Referring to Lemma \ref{lmae}, $M(T)$ has a pair of conjugate eigenvalues
$\lambda$, $\bar{\lambda}$ with $|\lambda|=e^{-cT/2}$. Thus $P(X)$ is a
contracting mapping. According to Lemma \ref{lamcc}, there is a $C^{1}$
diffeomorphism $\varphi$ which is near enough to the identity that $PX-X_{0}$
is conjugate equivalent to $M(T)$. There is an invertible constant matrix $C$
such that
\[
C^{-1}M(T)C=\left(
\begin{array}
[c]{cc}%
\lambda & 0\\
0 & \bar{\lambda}%
\end{array}
\right)  =D(\lambda),
\]
and we may suppose that
\begin{equation}
\frac{1}{2}\left\vert X-X_{0}\right\vert <\left\vert \varphi(X)-\varphi
(X_{0})\right\vert <2\left\vert X-X_{0}\right\vert \label{hbcm}%
\end{equation}
for $X-X_{0}$ small, since $\varphi$ is near the identity. Therefore, the
Liapunov exponent is given by
\begin{align*}
\mu_{x}  &  =\lim_{n\rightarrow\infty}\frac{1}{nT}\ln\left\vert P^{n}%
X-X_{0}\right\vert \\
&  =\lim_{n\rightarrow\infty}\frac{1}{nT}\ln\left\vert \varphi\circ
M(T)^{n}\circ\varphi^{-1}(X)-\varphi\circ M(T)^{n}\circ\varphi^{-1}%
(X_{0})\right\vert \\
&  =\lim_{n\rightarrow\infty}\frac{1}{nT}\ln\left\vert D(\lambda)^{n}\circ
C^{-1}\left[  \varphi^{-1}(X)-\varphi^{-1}(X_{0})\right]  \right\vert
=-\frac{c}{2}.
\end{align*}
The second equality follows from (\ref{hbcm}). Hence, the rate of decay of the
solution to the unique $T$-periodic solution is $c/2$, independently of the
initial value $X$.

\begin{rmk}
The above result shows that the Lyapunov exponent is invariant under a $C^{1}$
conjugate transformation. From the proof, the conclusion is still true if
$\varphi$ is a bi-Lipschitz mapping.
\end{rmk}

\subsection{\label{Prothm:a}Proof of Theorem \ref{thm:a}}

The existence and uniqueness of $T$-periodic solutions of (\ref{ii}) can be
obtained by the same argument in the proof of Theorem \ref{thmb}. But for the
stability of the periodic solution, the $C^{1}$ regularity of the solutions of
(\ref{ii}) with respect to the initial value is needed. The following result
is due to Lazer and McKenna.

\begin{lma}
\label{lamcb} Let $u(t,\xi,\eta)$ be the solution of the initial-value
problem
\[
\left\{
\begin{array}
[c]{l}%
u^{\prime\prime}+cu^{\prime}+a(t)u^{+}-b(t)u^{-}=h(t),\\
u(0)=\xi,\quad u^{\prime}(0)=\eta.
\end{array}
\right.
\]
Assume that $a(t), b(t) \in C_{T}$, and that $h\in C_{T}$ has a finite number
of zeros. Then, for $t\in[\,0,\bar{t}\,]$, the partial derivatives of $u$ and
$u^{\prime}$ with respect to $\xi,\eta$ exist and are continuous. Moreover,
if
\[
X(t)=\left\{
\begin{array}
[c]{cc}%
\frac{\partial u}{\partial\xi} & \frac{\partial u}{\partial\eta}\\
\frac{\partial u^{\prime}}{\partial\xi} & \frac{\partial u^{\prime}}%
{\partial\eta}%
\end{array}
\right\}  ,
\]
then
\[
X^{\prime}(t)=A(t)X(t),\quad X(0)=\operatorname{Id},
\]
where
\[
A(t)=\left\{
\begin{array}
[c]{cc}%
0 & 1\\
-p(t) & -c
\end{array}
\right\}
\]
and $p(t)=a(t)\chi_{+}+b(t)\chi_{-}$, where $\chi_{\pm}$ denotes the
characteristic function of the set $\{t\in[\,0,T\,],\,u(t)\gtrless0\}$.
\end{lma}

\begin{rmk}
Though Lemma \ref{lamcb} was originally stated only for both $a$ and $b$
constant, by carefully examining the proof in the appendix of \cite{LaMckk},
one may assert that Lemma \ref{lamcb} holds when $a(t)$ and $b(t)$ are continuous.

Now that $u(t)$ and $u^{\prime}(t)$ are known to be $C^{1}$ with respect to
the initial value, by reasoning along the same lines as in the proof of
Theorem \ref{thmb}, the stability and rate of decay of the periodic solutions
can be obtained.
\end{rmk}

\begin{ack}
We wish to thank the referee for his/her helpful comments and suggestions.
\end{ack}

\end{document}